\documentclass[preprint,12pt]{elsarticle}

\usepackage{graphicx}

\usepackage{amssymb}

\journal{Journal of Computational Physics}

\begin{document}

\begin{frontmatter}



\title{Beyond the mesh handling Maxwell's curl equations with an unconditionally leapfrog stable scheme}


\author[a]{Guido Ala}
\author[b]{Elisa Francomano\fnref{label1} }
\fntext[label1]{Corresponding author: Elisa Francomano\\ 
\textit{E-mail address}: elisa.francomano@unipa.it\\
\textit{Postal address}: Universit\`a degli Studi di Palermo, Viale delle Scienze, Ed. 6, 90128 Palermo, Italy\\
\textit{Phone number}: ++39 091 23842515}


\address[a]{Universit\`a degli Studi di Palermo, DEIM, Palermo, Italy}
\address[b]{Universit\`a degli Studi di Palermo, DICGIM, Palermo, Italy}

\begin{abstract}
Numerical solution of equations governing time domain simulations in computational electromagnetics, is usually based on grid methods in space and on explicit schemes for the time evolution. A predefined grid in the problem domain and a stability step size restriction must be accepted. Evidence is given that efforts need for overcoming these heavy constraints. 
Recently, the authors developed a meshless method to avoid the connective laws among the points scattered in the problem domain. Despite the good spatial properties, the numerical explicit integration used in the original formulation of the method provides,also in a meshless context, spatial and time discretization strictly interleaved and mutually conditioned. Afterwards, in this paper the stability condition is firstly addressed in a general way by allowing the time step increment get away from the minimum points spacing. Meanwhile, a formulation of the alternating direction implicit scheme for the evolution in time is combined with the meshless solver.  The formulation preserves the leapfrog marching on in time of the explicit integration scheme. The new method, not constrained by a gridding in space and unconditionally stable in time, is numerical assessed by different numerical simulations. Perfect matching layer technique is used in simulating open spatial problems; otherwise, a consistency restoring approach is introduced in treating truncation at finite boundary and irregular points distribution. Three case studies are investigated by achieving a satisfactory agreement comparing both numerical and analytical results.
\end{abstract}

\begin{keyword} ADI leapfrog method, meshless methods, smoothed particles electromagnetics.

\end{keyword}

\end{frontmatter}


\section{Introduction}
Meshless methods have recently emerged as numerical techniques for electromagnetic modelling
\cite{bib_1,bib_3,bib_4,bib_5,bib_13,bib_15,bib_18,bib_18b,bib_19,bib_20,bib_25,bib_33a,bib_36}. These methods do not require a predefined mesh, and use points
scattered in the problem domain avoiding the need of information on the position among them. One of
the most popular meshless method, the smoothed particle hydrodynamics (SPH) 
\cite{bib_22,bib_23,bib_24,bib_25,bib_28,bib_29,bib_30,bib_31} has
been recently reformulated by the authors for computational electromagnetics (CEM) problems. 
Maxwell's equations, which relate the electric and magnetic fields by means of a time dependent system of partial differential equations(PDEs) are considered \cite{bib_1,bib_2}. Due to the coupling nature of the electric and magnetic field
components, the points are in such a way that the electric points should be surrounded by those
of magnetic field and vice-versa.
The method named as SPEM is applied in approximating the space field variables of the time domain Maxwell’s curl equations by using a kernel representation working with a cluster of scattered nodes.
The differential operators applied to the field variables are transmitted to a kernel function which
characterizes the spatial discretization. Despite the good spatial properties of SPEM, the numerical
integration based on explicit method for time evolution, introduces a severe constraint. In fact, as for
the conventional grid based methods, SPEM has to satisfy the Courant-Friedrich-Levy stability
condition \cite{bib_11,bib_16,bib_34}. As well known, this condition becomes highly restrictive when the geometry
resolution forces the use of much finer space discretization than that is needed to solve for waves away
from material irregularities or interfaces. This problem is firstly investigated in the paper in SPEM
framework, in order to determine the largest allowed time increment.
Moreover, in technical literature, a great interest is addressed in finding numerical schemes which allow computational advantages together with unconditional stability \cite{bib_14,bib_32,bib_36,bib_37,bib_38}. Otherwise, an alternate implicit direction method is proposed for the time evolution of the electromagnetic field in the SPEM framework by preserving the leapfrog marching on in time of the explicit integration scheme. By considering a 2-stage splitting approach, the leapfrog implicit formulation is carried out  by suitably handling the formulas related to two adjacent time steps. By working with an implicit scheme, a sparse linear system is generated and it has to be solved at each time step. The system matrix involves the kernel derivatives, and it
can be assembled only once at the beginning of time integration process. Details on the matrices used in the problem are investigated.
The resulting algorithm is called as LAF-SPEM method. Accuracy degradation, due to a lack of particle consistency,  is also taken into account: a numerical corrective strategy, which allows to restore the consistency, without any modification of the smoothing kernel function and its derivatives, is employed  \cite{bib_1}.
Test problems are simulated to validate the proposed methodology. Simulations running with a
time step larger than that satisfying the CFL condition has been found to remain stable, so improving
the conventional SPEM algorithm.
The paper is organized as follows. In section 2 basic ideas on smoothed particle hydrodynamics method
are recalled in order to provide the original SPEM numerical scheme in solving time-domain Maxwell’s curl
equations. Investigations about the stability criterion are also included. In section 3 the LAF-SPEM
method is described and simulation results are reported and discussed in section 4 also by comparing
with standard FDTD computations and existing exact solution. A brief conclusion complete the paper.

\section{The original smoothed particle electromagnetics method}
The SPH scheme has been adopted by the authors to compute the time-domain Maxwell’s curl
equations by using particle approximation, so obtaining SPEM method \cite{bib_1,bib_2}. For the reader
convenience, a brief summary of SPEM fundamentals is reported in the following. Let consider the
standard Maxwell ' s curl equations in time domain, in a linear medium in the compact curl notation:
\begin{eqnarray}
\label{eq1}
\nabla \times \textbf{H}=\textbf{J}+\epsilon D_{t}\textbf{E}
\end{eqnarray}
\begin{center}
$\nabla \times \textbf{E}=-\mu D_{t} \textbf{H}$
\end{center}
where $\textbf{E}$ and  $\textbf{H}$ are the electric and magnetic vectors fields,  $\textbf{J}=\sigma \textbf{E}$ is the conduction current density vector, $\textit{t}$ as the time, $D_{t}$ the derivative operator, $\epsilon, \mu, \sigma $ the medium permittivity, permeability, and conductivity, respectively, which are supposed to be constant. 
The SPH method is based on an integral representation of the unknown function $\textit{u(x)}$:
\begin{eqnarray}
\label{eq2}
u(x)=\int_{\Omega}u(y)\delta(x-y)dy
\end{eqnarray}\\
where $\Omega \subset \Re^{d}$ is the problem domain and $\delta(\cdot)$ is the Dirac delta function.
The delta function is approximated by a smoothing kernel function $W(\cdot)$ depending on the spatial
coordinates and on the smoothing length parameter $\textit{h}$. The smoothing kernel function must verify the following properties:
\begin{enumerate}
\item[a.]  $\lim_{h \rightarrow 0} W(x-y,h)=\delta(x-y);$
\item[b.]  $W(x-y,h)>0 $  on  a sub-domain of $\Omega$ and $W(x-y,h)=0$ outside;
\item [c.] $\int_{\Omega}W(x-y,h)dy=1;$
\item [d.] $W(x-y,h)$ function value for a particle is monotonically decreasing\\  with increasing distance away from the particle.
\end{enumerate}

The discrete spatial formulation of equations (\ref{eq1}) has been obtained by using SPH particle
approximation: the integral representation is approximated over a set of particles placed in the problem
domain. The particles are the points in which the fields components are computed at each time step, by
using the information belonging to the neighbouring ones. The smoothing length \textit{h} defines the size of
the support domain of each particle \cite{bib_18a, bib_23,bib_24} . Particular attention has to be paid on the choice of \textit{h}. A too small \textit{h} can lead to an inaccurate solution because not enough particles may
be placed inside the influence area of a fixed particle $x_{i}$ ; a too large smoothing length, can smooth out
local properties and the particle approximation suffers too. From now on, for each referred fixed
particle $x_{i}$ , an own smoothing space length $h_{i}$ will be considered.
The SPH particle approximation of (\ref{eq2}) can be expressed as a linear combination of translates of a kernel basis function:

\begin{eqnarray}
\begin{array}{lll}
\label{eq3}
u(x_{i})=\sum\limits_{j=1}^{NP} u(x_{j})W(x_{i}-x_{j},h_{i})\Delta V_{j}\\
\end{array}
\end{eqnarray}
where \textit{NP} is the total number of particles in the influence domain of the fixed one $x_{i}$ and $\Delta V_{j}$ measures the surrounding media at the position of particle $x_{j}$. \textit{NP} is strictly related to the choice of the smoothing length $h_{i}$. In computing the two interlaved Maxwell's curl equations, particles for the electric field ($\textit{E}$-nodes) and magnetic field ($\textit{H}$-nodes) must be considered in the spatial domain. Due to the electric and magnetic coupling nature of Maxwell' s equations, spatial staggering of \textit{E}-nodes and \textit{H}-nodes is required. The \textit{E}-nodes and the \textit{H}-nodes are arranged so that each \textit{E}-node is surrounded by \textit{H}-nodes and each \textit{H}-node is surrounded by \textit{E}-nodes. The electric field update depends on the magnetic particles in the neighbourhood of the fixed electric particle and viceversa. In solving ($\ref{eq1}$), fields spatial derivatives have to be approximated by means of SPH. 
The differential operators on the fields components are transmitted to the smoothing kernel function: so the kernel approximation allows differential operations to be determined from the values of the function and the
derivatives of the kernel, rather than the derivatives of the function itself. However, this is really true
only when the influence domain of a fixed particle is placed entirely inside the problem domain; when
the influence domain overlaps with the geometry boundary the smoothing kernel \textit{W} is truncated, and a
resulting non-zero surface integral has to be opportunely treated in order to avoid the corruption of the numerical solution \cite{bib_23}.
In the following, for the sake of simplicity, and in agreement with the numerical results reported in
section 4, a lossless medium is considered. Thus, the two interleaved Maxwell's curl equations (\ref{eq2}) result in SPEM formulation as follows:
\begin{eqnarray}
\begin{array}{lll}
\label{eq4}
D_{t}H_{x}(r_{i}^H) &=-\frac{1}{\mu}\sum\limits_{j=1}^{NP_{H}}[E_{z}(r_{j}^E)D_{y}W(r_{i}^H-r_{j}^E,h_{i}^H)-\\
& \quad E_{y}(r_{j}^E)D_{z}W(r_{i}^H-r_{j}^E,h_{i}^H)] \Delta V_{j}
\end{array}
\end{eqnarray}
\begin{eqnarray*}
\begin{array}{lll}
D_{t}H_{y}(r_{i}^H) &=-\frac{1}{\mu}\sum\limits_{j=1}^{NP_{H}}[E_{x}(r_{j}^E)D_{z}W(r_{i}^H-r_{j}^E,h_{i}^H)-\\
& \quad E_{z}(r_{j}^E)D_{x}W(r_{i}^H-r_{j}^E,h_{i}^H)] \Delta V_{j}\\
\end{array}
\end{eqnarray*}
\begin{eqnarray*}
\begin{array}{lll}
D_{t}H_{z}(r_{i}^H) &=-\frac{1}{\mu}\sum\limits_{j=1}^{NP_{H}}[E_{y}(r_{j}^E)D_{x}W(r_{i}^H-r_{j}^E,h_{i}^H)-\\
& \quad E_{x}(r_{j}^E)D_{y}W(r_{i}^H-r_{j}^E,h_{i}^H)] \Delta V_{j}\\
\end{array}
\end{eqnarray*}

\begin{eqnarray}
\begin{array}{lll}
\label{eq5}
D_{t}E_{x}(r_{i}^E) &=\frac{1}{\epsilon}\sum\limits_{j=1}^{NP_{E}}[H_{z}(r_{j}^H)D_{y}W(r_{i}^E-r_{j}^H,h_{i}^E)-\\
& \quad  H_{y}(r_{j}^H)D_{z}W(r_{i}^E-r_{j}^H,h_{i}^E)] \Delta V_{j}\\
\end{array}
\end{eqnarray}

\begin{eqnarray*}
\begin{array}{lll}
D_{t}E_{y}(r_{i}^E)& =\frac{1}{\epsilon}\sum\limits_{j=1}^{NP_{E}}[H_{x}(r_{j}^H)D_{z}W(r_{i}^E-r_{j}^H,h_{i}^E)-\\
& \quad H_{z}(r_{j}^H)D_{x}W(r_{i}^E-r_{j}^H,h_{i}^E)] \Delta V_{j}\\
\end{array}
\end{eqnarray*}

\begin{eqnarray*}
\begin{array}{lll}
D_{t}E_{z}(r_{i}^E)& =\frac{1}{\epsilon}\sum\limits_{j=1}^{NP_{E}}[H_{y}(r_{j}^H)D_{x}W(r_{i}^E-r_{j}^H,h_{i}^E)-\\
& \quad H_{x}(r_{j}^H)D_{y}W(r_{i}^E-r_{j}^H,h_{i}^E)] \Delta V_{j}\\
\end{array}
\end{eqnarray*}
where $r_{i}^E$ and $r_{j}^H$ denote the positions of the \textit{i}-th and \textit{j}-th particle for the electric and magnetic field component, respectively; $NP_{E}$ ($NP_{H}$) is the number of neighbouring particles of the \textit{i}-th fixed electric
(magnetic) particle. The first derivative operator is indicated as $D_{x} ,D_{y} ,D_{z}$, \textit{x, y,z} as the spatial coordinates. A central finite difference algorithm is used to approximate time derivatives in ($\ref{eq4}$) and ($\ref{eq5}$).
A leapfrog scheme requiring only explicit updates, has been used for advancing the solution in
time: the electric field is computed at whole time step and the magnetic field at the half time step \cite{bib_1}. By indicating with $\Delta t$ the time step and with $ n=n\Delta t$ the generic step, a complete temporal step of (\ref{eq1}) evolves from \textit{n}-1/2 to \textit{n}+1, by involving both electric and magnetic field updates. The marching-on in time is actually from the step \textit{n}-1/2 to the step \textit{n}+1/2 when updating magnetic field while is from step \textit{n} to the step \textit{n}+1 when updating electric field:
\begin{eqnarray}
\begin{array}{lll}
\label{eq6}
H_{x}^{n+1/2}(r_{i}^H) &=H_{x}^{n-1/2}(r_{i}^H)-\frac{\Delta t }{\mu}\sum\limits_{j=1}^{NP_{H}}[E_{z}^{n}(r_{j}^E)D_{y}W(r_{i}^H-r_{j}^E,h_{i}^H)-\\
& \quad E_{y}^{n}(r_{j}^E)D_{z}W(r_{i}^H-r_{j}^E,h_{i}^H)] \Delta V_{j}
\end{array}
\end{eqnarray}

\begin{eqnarray*}
\begin{array}{lll}
 H_{y}^{n+1/2}(r_{i}^H) & =H_{y}^{n-1/2}(r_{i}^H)-\frac{\Delta t }{\mu}\sum\limits_{j=1}^{NP_{H}}[E_{x}^{n}(r_{j}^E)D_{z}W(r_{i}^H-r_{j}^E,h_{i}^H)-\\
& \quad E_{z}^{n}(r_{j}^E)D_{x}W(r_{i}^H-r_{j}^E,h_{i}^H)] \Delta V_{j}\\
\end{array}
\end{eqnarray*}
\begin{eqnarray*}
\begin{array}{lll}
H_{z}^{n+1/2}(r_{i}^H) &=H_{z}^{n-1/2}(r_{i}^H)-\frac{\Delta t }{\mu}\sum\limits_{j=1}^{NP_{H}}[E_{y}^{n}(r_{j}^E)D_{x}W(r_{i}^H-r_{j}^E,h_{i}^H)-\\
&\quad E_{x}^{n}(r_{j}^E)D_{y}W(r_{i}^H-r_{j}^E,h_{i}^H)] \Delta V_{j}\\
\end{array}
\end{eqnarray*}
\begin{eqnarray}
\begin{array}{lll}
\label{eq7}
E_{x}^{n+1}(r_{i}^E) &=E_{x}^{n}(r_{i}^E)+\frac{\Delta t }{\epsilon}\sum\limits_{j=1}^{NP_{E}}[H_{z}^{n+1/2}(r_{j}^H)D_{y}W(r_{i}^E-r_{j}^H,h_{i}^E)-\\
& \quad H_{y}^{n+1/2}(r_{j}^H)D_{z}W(r_{i}^E-r_{j}^H,h_{i}^E)] \Delta V_{j}\\
\end{array}
\end{eqnarray}
\begin{eqnarray*}
\begin{array}{lll}
E_{y}^{n+1}(r_{i}^E) &=E_{y}^{n}(r_{i}^E)+\frac{\Delta t }{\epsilon}\sum\limits_{j=1}^{NP_{E}}[H_{x}^{n+1/2}(r_{j}^H)D_{z}W(r_{i}^E-r_{j}^H,h_{i}^E)-\\
& \quad H_{z}^{n+1/2}(r_{j}^H)D_{x}W(r_{i}^E-r_{j}^H,h_{ i}^E)] \Delta V_{j}\\
\end{array}
\end{eqnarray*}

\begin{eqnarray*}
\begin{array}{lll}
E_{z}^{n+1}(r_{i}^E) &=E_{z}^{n}(r_{i}^E)+\frac{\Delta t }{\epsilon}\sum\limits_{j=1}^{NP_{E}}[H_{y}^{n+1/2}(r_{j}^H)D_{x}W(r_{i}^E-r_{j}^H,h_{i}^E)-\\
& \quad H_{x}^{n+1/2}(r_{j}^H)D_{y}W(r_{i}^E-r_{j}^H,h_{i}^E)] \Delta V_{j}\\
\end{array}
\end{eqnarray*}

As well-known, the explicit time integration scheme is subjected to the CFL condition, that for the EM problem is expressed as:

\begin{center}
$\Delta t \le  \frac{\Delta r}{c\sqrt{d}}$
\end{center}
where $\textit{c}$ is the propagation velocity of the physical phenomenon in the medium, \textit{d} is the geometric dimension of the problem, and $\Delta r$ is the minimum points spacing in the problem domain. In SPEM, the previous constraint depends on the smoothing length value; moreover, this last is strictly related to the spatial particles distances. Therefore it requires the time step to be proportional to the smallest spatial point resolution: 
 $\Delta t \le min_{i}(hi/c).$

The previous condition guarantees the numerical stability, but it can makes the SPEM impractical to use, especially when the distance between the two closest nodes becomes very small. In this section, some efforts are made to obtain a better maximum allowable time step. Namely, referring to the first case study reported in section 4, a propagation of a 2D $TM_{z}$ wave is considered.
In this case, the discrete formulation of relations (\ref{eq6}) and (\ref{eq7}) are compactly re-written, as follows:
\begin{eqnarray}
\begin{array}{lll}
\label{eq8}
\textbf{E}_{z}^{n+1}=\textbf{E}_{z}^{n}+\Delta t [\textbf{T}\textbf{H}_{y}^{n+1/2}-\textbf{U}\textbf{H}_{x}^{n+1/2}]\\
\end{array}
\end{eqnarray}

\begin{eqnarray}
\begin{array}{lll}
\label{eq9}
\textbf{H}_{x}^{n+1/2}=\textbf{H}_{x}^{n-1/2}-\Delta t \textbf{V}\textbf{E}_{z}^{n}\\
\end{array}
\end{eqnarray}

\begin{eqnarray}
\begin{array}{lll}
\label{eq10}
\textbf{H}_{y}^{n+1/2}=\textbf{H}_{y}^{n-1/2}+\Delta t \textbf{Z} \textbf{E}_{z}^{n}\\
\end{array}
\end{eqnarray}
\\
where $\textbf{E}_{z},\textbf{H}_{x},\textbf{H}_{y}$ are the vectors of the field components and $\textbf{T},\textbf{U},\textbf{V},\textbf{Z}$ as the matrices collecting the spatial field derivatives computed with SPEM formulation and the physical parameters of the problem domain. By substituting (\ref{eq9}) and (\ref{eq10}) in (\ref{eq8}), i.e.:
\begin{eqnarray}
\begin{array}{lll}
\label{eq11}
\textbf{E}_{z}^{n+1}=[\textbf{I}+\Delta t^2(\textbf{TZ}+\textbf{UV})]\textbf{E}_{z}^{n}+\Delta t [\textbf{TH}_{y}^{n-1/2}-\textbf{UH}_{x}^{n-1/2}]\\
\end{array}
\end{eqnarray}\\
the evolution in time is described by the relation $\textbf{X}^{n+1}=\textbf{S}\textbf{X}^{n}$, where:\\
\begin{center}
$\textbf{S}=
\left(
\begin{array}{ccc}
\textbf{I}+\Delta t^2(\textbf{TZ}+\textbf{UV})&-\Delta t\textbf{U}&\Delta t\textbf{T}\\
-\Delta t \textbf{V}&\textbf{I}&0\\
\Delta t \textbf{Z}&0&\textbf{I}
\end{array}
\right)$\\
\end{center} 
\begin{center}
$\textbf{X}^{n+1}=
\left(
\begin{array}{ccc}
\textbf{E}_{z}^{n+1}\\
\textbf{H}_{x}^{n+1/2}\\
\textbf{H}_{y}^{n+1/2}
\end{array}, 
\right),\,
\textbf{X}^{n}= 
\left(
\begin{array}{ccc}
\textbf{E}_{z}^{n}\\
\textbf{H}_{x}^{n-1/2}\\
\textbf{H}_{y}^{n-1/2}
\end{array}
\right)$\\
\end{center}
In order to $\textbf{X}^{n+1}$ be bounded, as n $\longrightarrow \infty$, the time increment must satisfy the following relation [27]:
\begin{eqnarray}
\begin{array}{c}
\label{eq12}
\Delta t< \frac{2}{\sqrt{\lambda_{MAX}(\textbf{S})}}\\
\end{array}
\end{eqnarray}\\
where $\lambda_{MAX}(\textbf{S})$ is the maximum  eigenvalue of the matrix $\textbf{S}$.
The (\ref{eq12}) gives an upper limit to the increment of the time step getting away from the minimum particles spacing. Furthermore, the restriction on $\Delta t$ limits the physical applications yet, so that in the following  an unconditionally stable time advance equations has been introduced. The original SPEM method has been modified by preserving the leapfrog evolution in time of the explicit time integration scheme.

\section{The leapfrog ADI-FDTD method }
In this section the leapfrog ADI-FDTD method \cite{bib_14} is generally described  and re-formulated for generating the  implicit leapfrog formulation of SPEM. In the ADI-FDTD method, the temporal stepping from \textit{n}-1 to \textit{n} step is obtained by splitting the integer time step into two half-time steps. By indicating with \textit{E} and \textit{H} the fields components with respect to a generic coordinate in the \textit{x,y,z} geometric system, the finite difference time domain scheme involves the following general relations, in which the actual field value is formally expressed as a generic function of the fields at previous time steps:

\begin{eqnarray}
\begin{array}{lll}
\label{eq13}
\textit{E}^{n-1/2}=f(\textit{E}^{n-1},\textit{H}^{n-1/2}, \textit{H}^{n-1} )
\end{array}
\end{eqnarray}

\begin{eqnarray}
\begin{array}{lll}
\label{eq14}
\textit{H}^{n-1/2}=g(\textit{H}^{n-1},\textit{E}^{n-1/2},\textit{E}^{n-1} )
\end{array}
\end{eqnarray}

\begin{eqnarray}
\begin{array}{lll}
\label{eq15}
\textit{E}^{n}=f(\textit{E}^{n-1/2},\textit{H}^{n}, \textit{H}^{n-1/2})
\end{array}
\end{eqnarray}

\begin{eqnarray}
\begin{array}{lll}
\label{eq16}
\textit{H}^{n}=g(\textit{H}^{n-1/2},\textit{E}^{n}, \textit{E}^{n-1/2})
\end{array}
\end{eqnarray}\\

For the marching-on in time of the electric field components, equation (\ref{eq16}) is substituted in (\ref{eq15}) and (\ref{eq16}) is substituted in (\ref{eq15}), so obtaining:
\begin{eqnarray}
\begin{array}{lll}
\label{eq17}
\textit{E}^{n-1/2}=\phi(\textit{E}^{n-1},\textit{H}^{n-1} )
\end{array}
\end{eqnarray}

\begin{eqnarray}
\begin{array}{lll}
\label{eq18}
\textit{E}^{n}=\phi(\textit{E}^{n-1/2},\textit{H}^{n-1/2} )
\end{array}
\end{eqnarray}

In (\ref{eq17}) and (\ref{eq18}) first and second order spatial field derivatives hold and tridiagonal linear systems have to be solved to generate $\textit{E}^{n-1/2}$ and $\textit{}E^{n}$. 
An ADI-FDTD time step is so described by means of relations (\ref{eq17}) and (\ref{eq14}), and relation (\ref{eq18}) and 
(\ref{eq16}), respectively.
By considering the temporal step for the fields evolution from step $\textit{n}-1/2$ to $\textit{n}$ and obtaining $\textit{H}^{n}-1/2$ from (\ref{eq16}), the following relation holds: 
\begin{eqnarray}
\begin{array}{lll}
\label{eq19}
\textit{H}^{n-1/2}=$\~{g}$(\textit{H}^{n},\textit{E}^{n},\textit{E}^{n-1/2} ).
\end{array}
\end{eqnarray}

By substituting (\ref{eq19}) into (\ref{eq15}) the electric field component is generated by means of: \\

\begin{eqnarray}
\begin{array}{lll}
\label{eq20}
\textit{E}^{n}=\psi(\textit{E}^{n-1/2},\textit{H}^{n})
\end{array}
\end{eqnarray}

Let now consider equation (\ref{eq17}) for the marching-on in time of the electric field from \textit{n} to $\textit{n}+1/2$:\\

\begin{eqnarray*}
\begin{array}{lll}
\textit{E}^{n+1/2}=\phi(\textit{E}^{n},\textit{H}^{n})
\end{array}
\end{eqnarray*}

and by using (\ref{eq20}), the previous relation can be re-written as:\\

\begin{eqnarray}
\begin{array}{lll}
\label{eq21}
\tilde{\textit{E}}^{n+1/2}=\tilde{\phi}( \textit{E}^{n-1/2},\textit{H}^{n} )
\end{array}
\end{eqnarray}\\
The leapfrog equation describing the marching-on in time of the electric field component is so obtained.
Let now consider relation (\ref{eq15}) for the half time step from \textit{n}+1/2 to \textit{n}+1:

\begin{eqnarray*}
\begin{array}{lll}
\textit{E}^{n+1}=f(\textit{E}^{n+1/2},\textit{H}^{n+1},\textit{H}^{n+1/2} )
\end{array}
\end{eqnarray*}\\

By substituting it into (\ref{eq16}), the following relation holds
\begin{eqnarray}
\begin{array}{lll}
\label{eq22}
\tilde{\textit{H}}^{n+1}=g( \textit{H}^{n+1/2},\textit{E}^{n+1}, \textit{E}^{n+1/2})=\gamma(\textit{H}^{n+1/2},\textit{E}^{n+1/2})
\end{array}
\end{eqnarray}\\ 
For the half time step from \textit{n} to \textit{n}+1/2 relation (\ref{eq13}) is considered:\\
\begin{eqnarray*}
\begin{array}{lll}
\textit{E}^{n}=\tilde f (\textit{E}^{n+1/2},\textit{H}^{n+1/2},\textit{H}^{n} )
\end{array}
\end{eqnarray*}\\

By substituting this last into (\ref{eq14}). the $\textit{H}^{n+1/2}$ component holds:\\
\begin{eqnarray}
\begin{array}{lll}
\label{eq23}
\textit{H}^{n+1/2}=g(\textit{H}^{n},\textit{E}^{n+1/2},\tilde{f}(\textit{E}^{n+1/2},\textit{H}^{n+1/2},\textit{H}^{n} ))
\end{array}
\end{eqnarray}\\ 
This last relation is subsituted in (\ref{eq23}) so obtaining the following:\\
\begin{eqnarray}
\begin{array}{lll}
\label{eq24}
\textit{H}^{n+1}=\tilde{\gamma}( \textit{H}^{n},\textit{E}^{n+1/2}).
\end{array}
\end{eqnarray}\\ 
Equations (\ref{eq21}) and (\ref{eq23}) give rise to LAF scheme for marching-on in time of the generic magnetic field component \textit{H}. In (\ref{eq21}) and (\ref{eq23}) only first and second spatial field derivatives are involved. The electric and magnetic fields are staggered in time and updated in one full time step, iteratively.
\subsection{The leapfrog ADI-FDTD method in SPEM}
The complete formulation of Maxwell's curl equations derived from eqs (\ref{eq21}) and (\ref{eq24}) is reported in the following by approximating the spatial first and second derivatives by means of SPH particle approximation:
\begin{eqnarray}
\begin{array}{lll}
\label{eq25}
\epsilon{E}_{x}^{n+1/2}(r_{i}^{E})-(\frac{\Delta t}{2})^{2}\frac{1}{\mu}\sum\limits_{j=1}^{NP_{E}}E_{x}^{n+1/2}(r_{j}^E)D_{y}^{2}W(p_{j}^E,h_{i}^{E})\Delta V_{j}=\\
\epsilon{E}_{x}^{n-1/2}(r_{i}^{E})-(\frac{\Delta t}{2})^{2}\frac{1}{\mu}
\sum\limits_{j=1}^{NP_{H}}E_{x}^{n-1/2}(r_{j}^E)D_{y}^{2}W(p_{j}^E,h_{i}^{E})\Delta V_{j}+\\
\Delta t \sum\limits_{j=1}^{NP_{H}}[H_{z}^{n}(r_{j}^H)D_{y}W(p_{j}^H,h_{i}^{E})-
H_{y}^{n}(r_{j}^H)D_{z}W(p_{j}^H,h_{i}^{E})]\Delta V_{j}
\end{array}
\end{eqnarray}

\begin{eqnarray}
\begin{array}{lll}
\label{eq26}
\epsilon{E}_{y}^{n+1/2}(r_{i}^{E})-(\frac{\Delta t}{2})^{2}\frac{1}{\mu}\sum\limits_{j=1}^{NP_{E}}E_{y}^{n+1/2}(r_{j}^E)D_{z}^{2}W(p_{j}^E,h_{i}^{E})\Delta V_{j}=\\
\epsilon{E}_{y}^{n-1/2}(r_{i}^{E})-(\frac{\Delta t}{2})^{2}\frac{1}{\mu}
\sum\limits_{j=1}^{NP_{E}}E_{y}^{n-1/2}(r_{j}^E)D_{z}^{2}W(p_{j}^E,h_{i}^{E})\Delta V_{j}+\\
\Delta t \sum\limits_{j=1}^{NP_{H}}[H_{x}^{n}(r_{j}^H)D_{z}W(p_{j}^H,h_{i}^{E})-
H_{z}^{n}(r_{j}^H)D_{x}W(p_{j}^H,h_{i}^{E})]\Delta V_{j}
\end{array}
\end{eqnarray}

\begin{eqnarray}
\begin{array}{lll}
\label{eq27}
\epsilon{E}_{z}^{n+1/2}(r_{i}^{E})-(\frac{\Delta t}{2})^{2}\frac{1}{\mu}\sum\limits_{j=1}^{NP_{E}}E_{z}^{n+1/2}(r_{j}^E)D_{x}^{2}W(p_{j}^E,h_{i}^{E})\Delta V_{j}=\\
\epsilon {E}_{z}^{n-1/2}(r_{i}^{E})-(\frac{\Delta t}{2})^{2}\frac{1}{\mu}
\sum\limits_{j=1}^{NP_{E}}E_{z}^{n-1/2}(r_{j}^E)D_{x}^{2}W(p_{j}^E,h_{i}^{E})\Delta V_{j}+\\
\Delta t \sum\limits_{j=1}^{NP_{H}}[H_{y}^{n}(r_{j}^H)D_{x}W(p_{j}^H,h_{i}^{E})-
H_{x}^{n}(r_{j}^H)D_{y}W(p_{j}^H,h_{i}^{E})]\Delta V_{j}
\end{array}
\end{eqnarray}

\begin{eqnarray}
\begin{array}{lll}
\label{eq28}
\epsilon{H}_{x}^{n+1}(r_{i}^{H})-(\frac{\Delta t}{2})^{2}\frac{1}{\mu}\sum\limits_{j=1}^{NP_{H}}H_{x}^{n+1}(r_{j}^H)D_{y}^{2}W(q_{j}^H,h_{i}^{H})\Delta V_{j}=\\
\epsilon {H}_{x}^{n}(r_{i}^{H})-(\frac{\Delta t}{2})^{2}\frac{1}{\mu}
\sum\limits_{j=1}^{NP_{H}}H_{x}^{n}(r_{j}^H)D_{y}^{2}W(q_{j}^H,h_{i}^{H})\Delta V_{j}+\\
\frac{\Delta t}{\mu} \sum\limits_{j=1}^{NP_{E}}[E_{y}^{n+1/2}(r_{j}^E)D_{z}W(q_{j}^E,h_{i}^{H})- E_{z}^{n+1/2}(r_{j}^E)D_{y}W(q_{j}^E,h_{i}^{H})]\Delta V_{j}
\end{array}
\end{eqnarray}

\begin{eqnarray}
\begin{array}{lll}
\label{eq29}
\epsilon{H}_{y}^{n+1}(r_{i}^{H})-(\frac{\Delta t}{2})^{2}\frac{1}{\mu}\sum\limits_{j=1}^{NP_{H}}H_{y}^{n+1}(r_{j}^H)D_{z}^{2}W(q_{j}^H,h_{i}^{H})\Delta V_{j}=\\
\epsilon {H}_{y}^{n}(r_{i}^{H})-(\frac{\Delta t}{2})^{2}\frac{1}{\mu}
\sum\limits_{j=1}^{NP_{H}}H_{y}^{n}(r_{j}^H)D_{z}^{2}W(q_{j}^H,h_{i}^{H})\Delta V_{j}+\\
\frac{\Delta t}{\mu} \sum\limits_{j=1}^{NP_{E}}[E_{z}^{n+1/2}(r_{j}^E)D_{x}W(q_{j}^E,h_{i}^{H})-
E_{x}^{n+1/2}(r_{j}^E)D_{z}W(q_{j}^E,h_{i}^{H}) ]\Delta V_{j}
\end{array}
\end{eqnarray} 

\begin{eqnarray}
\begin{array}{lll}
\epsilon{H}_{z}^{n+1}(r_{i}^{H})-(\frac{\Delta t}{2})^{2}\frac{1}{\mu}\sum\limits_{j=1}^{NP_{H}}H_{z}^{n+1}(r_{j}^H)D_{x}^{2}W(q_{j}^H,h_{i}^{H})\Delta V_{j}=\\
\epsilon {H}_{z}^{n}(r_{i}^{H})-(\frac{\Delta t}{2})^{2}\frac{1}{\mu}
\sum\limits_{j=1}^{NP_{H}}H_{z}^{n}(r_{j}^H)D_{x}^{2}W(q_{j}^H,h_{i}^{H})\Delta V_{j}+\\
\frac{\Delta t}{\mu} \sum\limits_{j=1}^{NP_{E}}[E_{x}^{n+1/2}(r_{j}^E)D_{y}W(q_{j}^E,h_{i}^{H})-
E_{y}^{n+1/2}(r_{j}^E)D_{x}W(q_{j}^E,h_{i}^{H})]\Delta V_{j}
\end{array}
\end{eqnarray}
where $p_{j}^E=r_{i}^E-r_{j}^E$, $p_{j}^H=r_{i}^E-r_{j}^H$, $q_{j}^H=r_{i}^H-r_{j}^H$, $q_{j}^E=r_{i}^H-r_{j}^E$. \\
The equations (25)-(30) can be expressed in a compact matrix notation as follows:
\begin{eqnarray}
\begin{array}{lll}
\label{eq31}
\textbf{A}\textbf{E}^{n+1/2}=\textbf{A}\textbf{E}^{n-1/2}+\textbf{B}\textbf{H}^{n}\\
\textbf{C}\textbf{H}^{n+1}=\textbf{C}\textbf{H}^{n}+\textbf{F}\textbf{E}^{n+1/2}
\end{array}
\end{eqnarray}\\ 
where:

\begin{center}
$\textbf{A}=\left(
\begin{array}{ccc}
\textbf{A}_{xx} &\textbf{0} &\textbf{0}\\
\textbf{0}  &\textbf{A}_{yy} &\textbf{0}\\
\textbf{0}  &\textbf{0}  &\textbf{A}_{zz}\\
\end{array}
\right),\,
\textbf{B}=\left(
\begin{array}{ccc}
\textbf{0}  &-\textbf{B}_{xy} &\textbf{B}_{xz}\\
\textbf{B}_{yx}  &\textbf{0}  &-\textbf{B}_{yz}\\
-\textbf{B}_{zx}  &\textbf{B}_{zy}  &\textbf{0}\\
\end{array}
\right)$
\end{center}
\begin{center}
$\textbf{C}=\left(
\begin{array}{ccc}
\textbf{C}_{xx} &\textbf{0} &\textbf{0}\\
\textbf{0}  &\textbf{C}_{yy} &\textbf{0}\\
\textbf{0}  &\textbf{0}  &\textbf{C}_{zz}\\
\end{array}
\right),\,
\textbf{F}=\left(
\begin{array}{ccc}
\textbf{0}  &-\textbf{F}_{xy} &\textbf{F}_{xz}\\
\textbf{F}_{yx}  &\textbf{0}  &-\textbf{F}_{yz}\\
-\textbf{F}_{zx}  &\textbf{F}_{zy}  &\textbf{0}\\
\end{array}
\right)$
\end{center}
and the non zero entries are: 

\[\begin{array}{*{20}{c}}
{k = x,y,z;}\\
{\beta  = y,z,x}
\end{array}\left\{ {\begin{array}{*{20}{c}}
{{a_{kk}}\left( {i,j} \right) = \varepsilon  - {{\left( {\frac{{\Delta t}}{2}} \right)}^2}\frac{1}{\mu }D_\beta ^2W\left( {p_j^E,h_i^E} \right)\Delta {V_j}}&{}\\
{{b_{kl}}\left( {i,j} \right) = {b_{lk}}\left( {i,j} \right) = \Delta t{D_\eta }W\left( {p_j^H,h_i^E} \right)\Delta {V_j}}&\begin{array}{l}
l = x,y,z; l \ne k\\
\eta  = z,y,x; \eta \ne k
\end{array}
\end{array}} \right.\quad \begin{array}{*{20}{c}}\\
\end{array}\]

\[\begin{array}{*{20}{l}}
{k = x,y,z;}\\
{\beta  = y,z,x}
\end{array}\left\{ {\begin{array}{*{20}{c}}
{{c_{kk}}\left( {i,j} \right) = \varepsilon  - {{\left( {\frac{{\Delta t}}{2}} \right)}^2}\frac{1}{\mu }D_\beta ^2W\left( {q_j^H,h_i^H} \right)\Delta {V_j}}&{}\\
{{f_{kl}}\left( {i,j} \right) = {f_{lk}}\left( {i,j} \right) = \frac{{\Delta t}}{\mu }{D_\eta }W\left( {q_j^E,h_i^H} \right)\Delta {V_j}}&\begin{array}{l}
l = x,y,z;{\rm{  }}l \ne k\\
\eta  = z,y,x; \eta  \ne k
\end{array}\\
\end{array}} \right.\quad \begin{array}{*{20}{c}}
\end{array}\]


The matrices employed in \textbf{A},\textbf{B},\textbf{C},\textbf{F} are sparse block matrices and the sparsity depends on the influence domain of each particle. Each row of the matrices refers to a fixed particle, center of a kernel function, and the non zero elements correspond to the particles surrounding the fixed one. The matrix $\textbf{A}_{xx}$  ($\textbf{A}_{yy},\textbf{A}_{zz}$) is generated from the smoothing kernel functions regarding the \textit{E}-particles , i.e. each kernel is centered on  $r_{i}^E$ and has $r_{j}^E$ as neighbourings. In a similar fashion the matrix $\textbf{C}_{xx}$  ($\textbf{C}_{yy},\textbf{C}_{zz}$) is related to the smoothing kernel functions centered on  $r_{i}^H$ and has $r_{j}^H$ as neighbourings. In all these matrices the second derivatives are involved.
The matrices $\textbf{B}_{kl}$ ($ k=x,y,z; l=x,y,z, l \neq k$) regard smoothing kernel centered on  $r_{i}^E$ but involving  $r_{j}^H$ as neighbourings. Analougously,  $\textbf{F}_{kl}$ ($ k=x,y,z; l=x,y,z, l \neq k $)  regards smoothing kernel fixed on  $r_{i}^H$ but involving  $r_{j}^E$ as neighbourings. In all these last matrices the first derivatives are involved.
When kernel vanish,  the limit configuration of the explicit formulation in time is obtained. In the following, some matrices configurations are reported. Equally distributed particles  are considered  generating a structured block matrices (Figs.1,2). In this case the matrices are banded and the bandwidth is with the width of the local support of the kernel function. The smoothing length  \textit{h}=$\alpha\Delta r $ must be opportunely chosen to achieve satisfactory approximation field values, giving rise to a sparse and solvable system.
The approach in which the support size remain fixed for increasingly denser set of data is very inefficient with associate system matrices increasingly denser. In the simulations, the smoothing length is scaled so that peaked basis functions are with densely spaced particles and flat basis functions are with coarsely spaced particles. \\In Fig.3 the error in the computation is reported by considering 100 equally distributed points over the whole problem domain. For the electric field in the $TM_{z}$ case, the error is measured in $L_{2}$ in comparison with the finite difference time domain solver. The best behaviour is with $\alpha$=0.7075. In tab.1 results are reported by considering the optimal smoothing length found, to give an accurate solution for this problem.  The sparsity of the matrices increases by increasing the problem size.

\begin{center}
\begin{tabular}{c|c}
{\textbf{N}\/} &{\textbf{Sparsity}\/}\\ \hline
100 & 81.86 \%  \\ 
200 & 91.412 \%  \\ 
625 & 96.45 \%   \\ 
900 & 96.45 \% \\ 
2500 & 97.48\% \\
3600 & 99.33\% \ 
\end{tabular}
\end{center}
\textbf{Table 1}. Results with regular particle distribution by increasing the particles density N over the problem domain in the computation of the matrix $\textbf{A}_{xx}$ \\

\begin{figure}
\label{fig:fig1}
\center
\includegraphics[width=10cm]{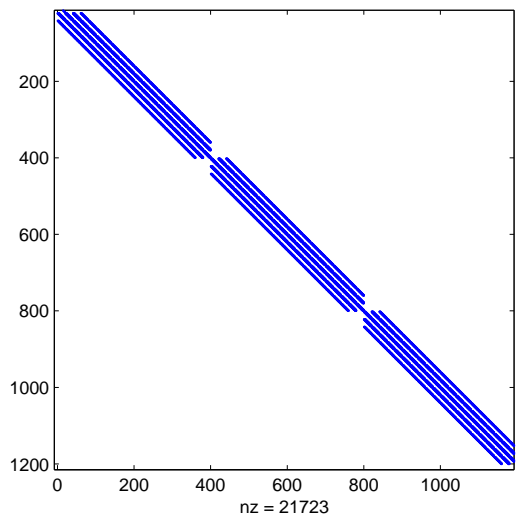}
\caption{Profile of the matrix $\textbf{A}$ for a regular particle distribution}
\end{figure}

\begin{figure}
\label{fig:fig2}
\center
\includegraphics[width=10cm]{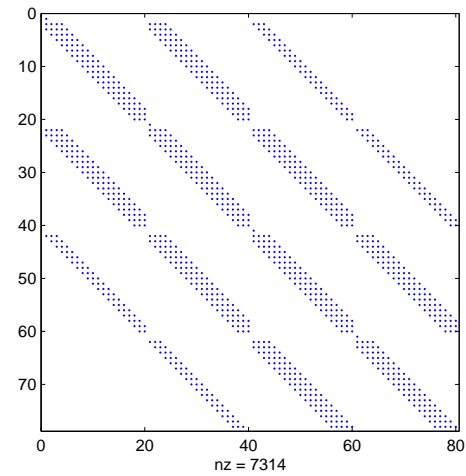}
\caption{Zoom of the skeleton of the matrix $\textbf{A}_{xx}$ in $\textbf{A}$ of figure 1}
\end{figure}

\begin{figure}[ht]
\label{fig:fig3}
\center
\includegraphics[width=10cm]{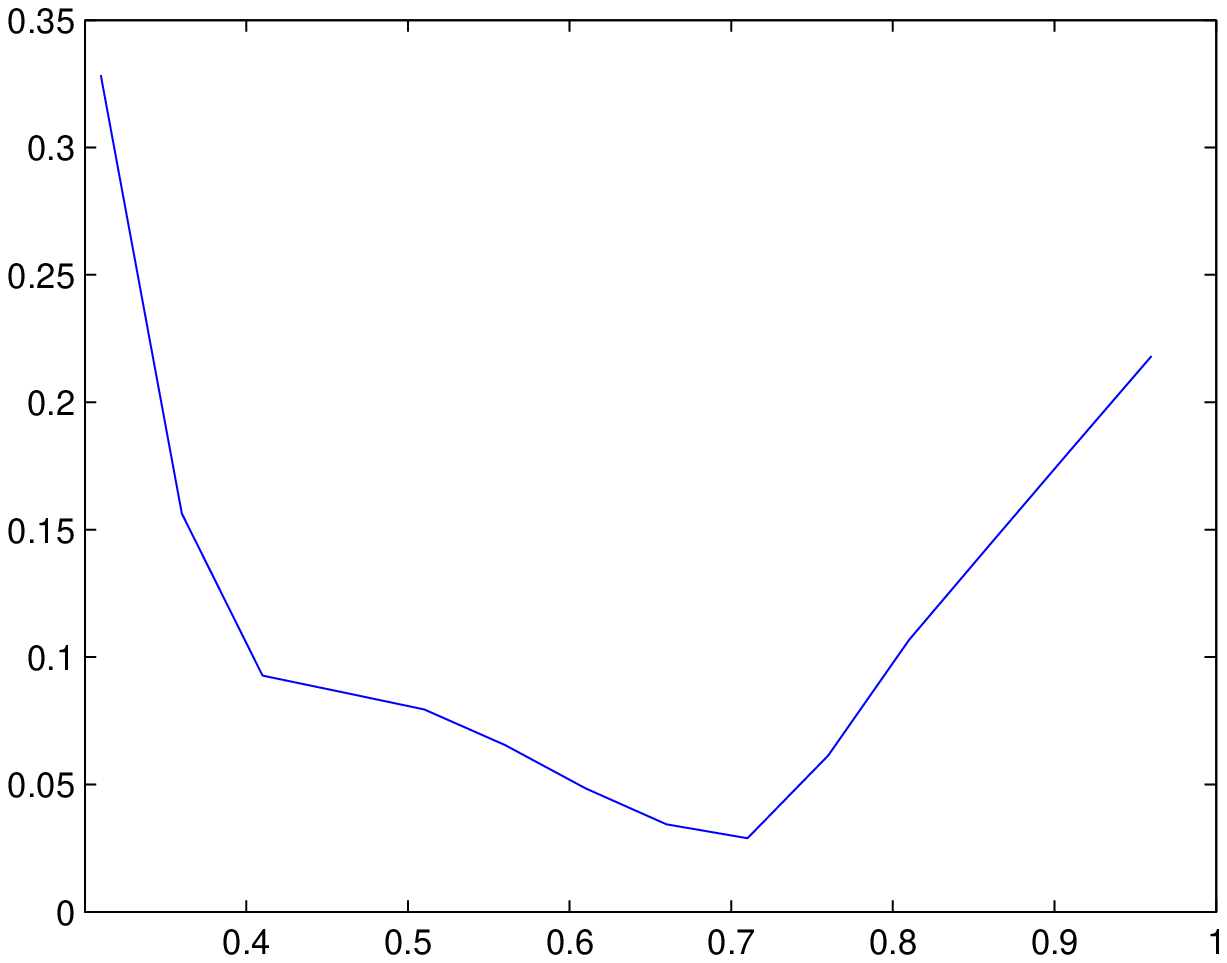}
\caption{Error behaviour by varying \textit{h} for 100 equally distributed particles over the problem domain.}
\end{figure}

In Figs.4,5 the block matrix $\textbf{A}_{xx}$ is depicted for irregular particle distributions, by varying the number of particle in the problem domain. Fig.4 is with 81 particles and a sparsity of about 80 \% is generated. Fig.5 refers to 625 particles with sparse matrix of about 96 \% of zero elements.

\begin{figure}
\label{fig:fig4}
\center
\includegraphics[width=10cm]{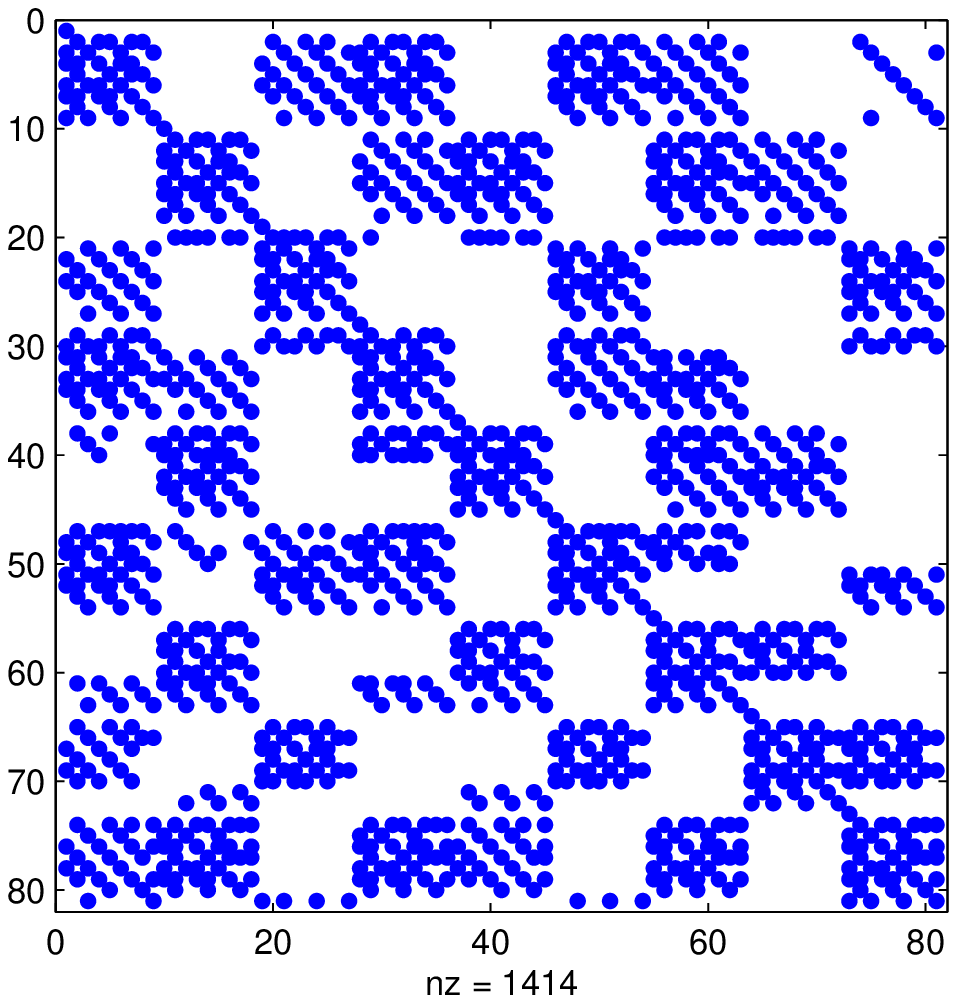}
\caption {Profile of the matrix $\textbf{A}_{xx}$ for an irregular 9x9 particle distribution}
\end{figure}
\begin{figure}
\label{fig:fig5}
\center
\includegraphics[width=10cm]{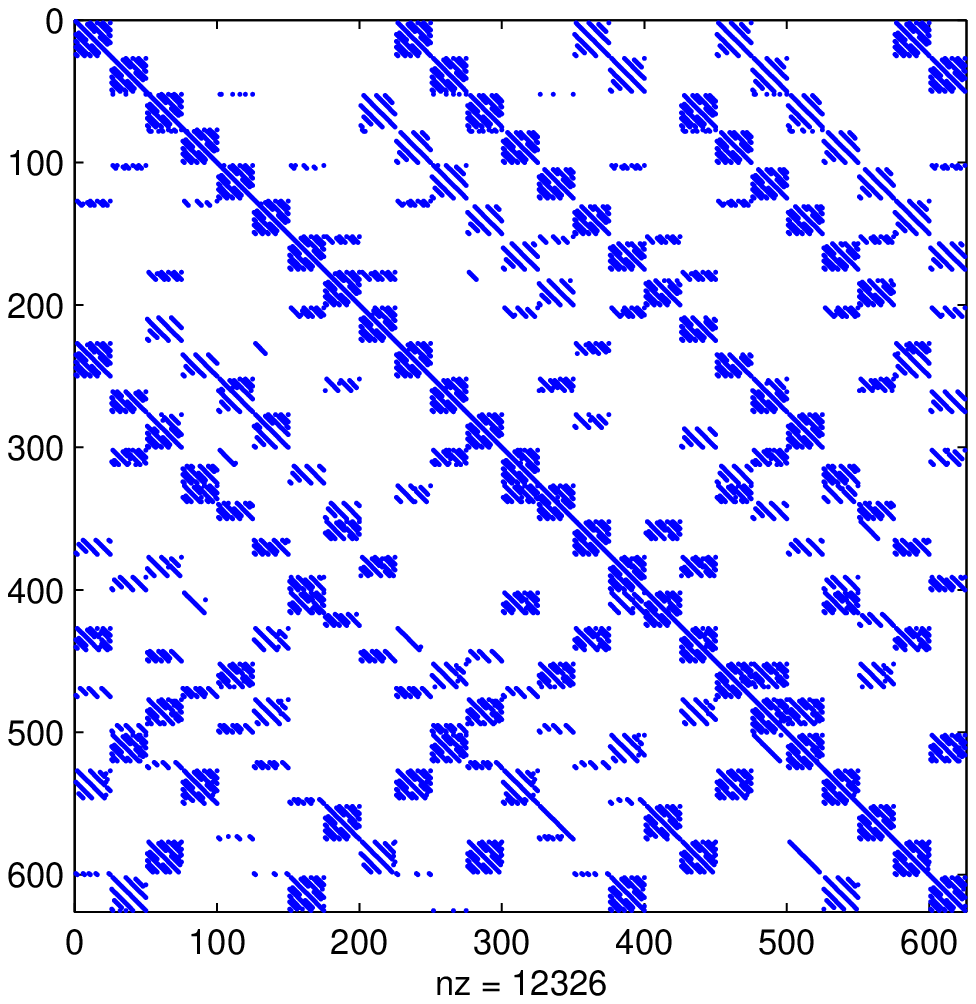}
\caption {Profile of the matrix $\textbf{A}_{xx}$ for an irregular 25x25 particle distribution
}
\end{figure}

The kernel derivatives in computing the derivatives of the field components can be performed in a pre-processing stage, and all the matrices employed in the computation can be assembled out of the temporal loop. Moreover, by considering the smoothing length independent from the evolution in time, the structure of the matrix not depends on the temporal step yet.

\section{Numerical validation}
In this section, some different case studies have been performed and results are discussed to assess the proposed meshless implicit leapfrog numerical method. At first, a transverse electric $(TM_{z})$) mode in air is firstly taken into account. In this case, equations are reduced as follows:

\begin{eqnarray}
\begin{array}{ccc}
\label{eq32}
D_{t} E_{z}=\frac{1}{\epsilon} (D_{x} H_{y}-D_{y} H_{x}), \,  D_{t} H_{y}=\frac{1}{\mu} (D_{x} E_{z}), \,   D_{t} H_{x}=\frac{1}{\mu} (D_{y} E_{z})
\end{array}
\end{eqnarray}

A regular particles distribution is built up, such as in a FDTD spatial grid. For problem with regular distribution, nodes are placed as in the point-matching procedure; for the general case, a Voronoi decomposition \cite{bib_34}or others technique can be adopted. In particular, in the paper when an irregular particles distribution is set, in order to avoid particle inconsistency which can lead to a loss of accuracy the consistency restoring adopted in \cite{bib_1,bib_2}, is used. The FDTD simulation is considered for comparison. The perfect matching layer (PML) technique is used in simulating the open spatial domain \cite{bib_16, bib_34}. In this way, also a good treatment of the surface integral truncation problem, as reported in section 2, is carried out for LAF-SPEM computation. The actual time step, satisfies the CFL condition, so as in the FDTD approach, i.e.:

\begin{eqnarray}
\begin{array}{lll}
\label{eq33}
\Delta t \equiv \Delta _{CFL}=\frac{\Delta x}{2c}\\
\end{array}
\end{eqnarray}
where $\Delta x=\Delta y= 1cm \,$ is the spatial step.
In figure 6 FDTD and LAF-SPEM spatial profiles along y direction, ten cells distant from the middle in the x direction in a 2D domain, 100x100 number of cells are reported; a time varying shaped source, $e^{\frac{(20-t)^{2}}{72}}$ , placed at the middle point of the domain is considered in simulating a $TM_{z}$ field in air. The comparison shows a satisfactory agreement.\\

\begin{figure}[ht]
\label{fig:fig6}
\center
\includegraphics[width=12cm]{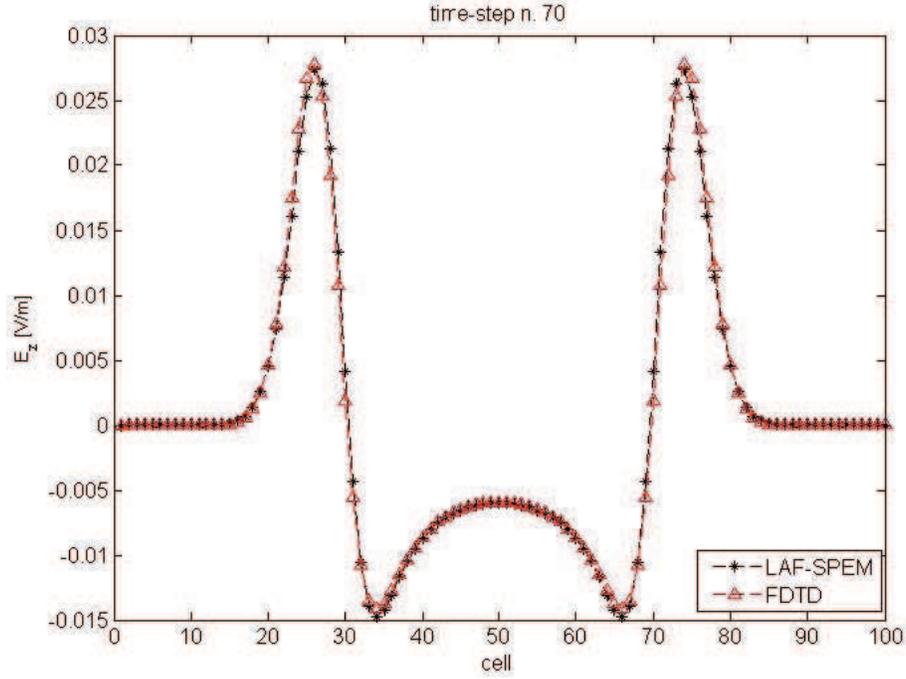}
\caption{Comparison among FDTD and LAF-SPEM results in a 2-D domain 100x100 cells and particles: the time dependent source is placed at the middle point and a $TM_{z}$ field is simulated in air. $\textbf{E}_{z}$ field (V/m) spatial profiles at time step 70, along y direction ten cells distant from the middle in the x direction}
\end{figure}

As a further validation, an axial symmetric cylindrical domain in air is considered, with $r_{0}=0.2m$ ($\epsilon_{0}, \mu_{0}$ as constitutive parameters), with the following boundary and initial conditions:
\begin{eqnarray}
\begin{array}{lll}
\label{eq34}
r=\sqrt{x^{2}+y^{2}}, r_{0}=0.2m,  \, 0 \le r \le r_{0}
\end{array}
\end{eqnarray}\\ 
\begin{eqnarray}
\begin{array}{lll}
\label{eq35}
E_{z}(x,y,0)=1-\frac{r^{2}}{r_{0}^{2}}, \,E_{z}(x_{0},y_{0},t)=0, \, \frac{\delta E_{z}(x,y,0)}{\delta t}=0
\end{array}
\end{eqnarray}\\ 

The domain is discretized with 10x10 irregular distributed particles obtained by randomly positioning the particles near the original regular position in the cylindrical domain. The actual time step is set equal to  $4 \Delta t_{CFL}$. As shown in figure 7 the simulation remains stable also for a final time corresponding to 175 time steps, and a satisfactory approximation has been achieved. 
\begin{figure}[ht]
\label{fig:fig7}
\center
\includegraphics[width=12cm]{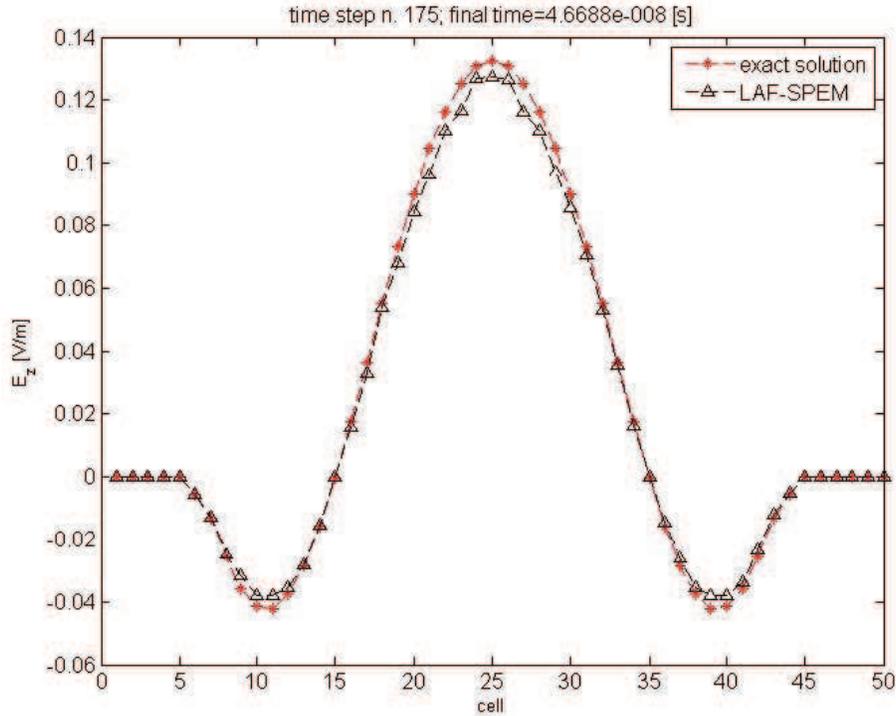}
\caption{Cylindrical domain with $r_{0}=0.2 m$ bounded by a perfect electric conductor: comparison between LAF-SPEM results and the exact solution at time step 175. The LAF-SPEM time step is 4 times that satisfying the CFL condition.}
\end{figure}

In order to better assess the validity of the proposed approach, a more realistic example has been carried out. A sectorial (2-D) Perfectly Electric Conducting (PEC) horn antenna excited by a sinusoidal voltage in a transverse magnetic $(TM_{z})$ computational domain has been simulated and compared with FDTD results [16] 
\cite{bib_16}. The computational domain is truncated with a PML absorbing boundary condition. The horn is modelled by a perfect electric conductor (PEC) material. The spatial cell size is equal to 0.0025 m, the time step is equal to 4.23 ps, the excitation frequency is equal to 9.84 GHz. 

\begin{figure}[ht]
\label{fig:fig8}
\center
\includegraphics[width=12cm]{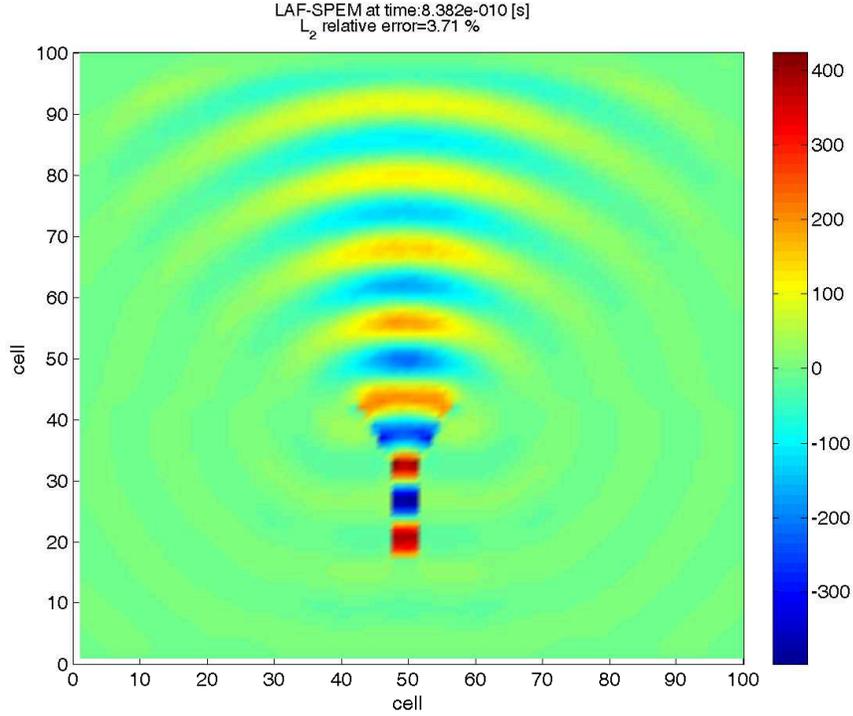}
\caption{ Sectorial (2-D) PEC horn antenna excited by a sinusoidal voltage in a 2-D domain with 100x100 regular particles distribution. Colour map of the Ez field (V/m) at time step 198 obtained with LAF-SPEM - relation (\ref{eq33})holds. The reported relative error is related to FDTD results.}
\end{figure}

In figure 8, LAF-SPEM $\textbf{E}_{y}$ field results, in a 2-D domain with 100x100 regular particles distribution, are reported, at time step 198. The excitation section is also reported in the figure. The actual time step, $\Delta t$ , satisfies the CFL condition (\ref{eq33}). In comparison with FDTD simulation, a good agreement has been obtained. For a better comparison with FDTD results, in figure 9 FDTD and LAF-SPEM $\textit{E}_{y}$ field spatial profiles along the direction containing the longitudinal axis of the horn antenna at the same time step of figure 6, are reported.

\begin{figure}[ht]
\label{fig:fig9}
\center
\includegraphics[width=12cm]{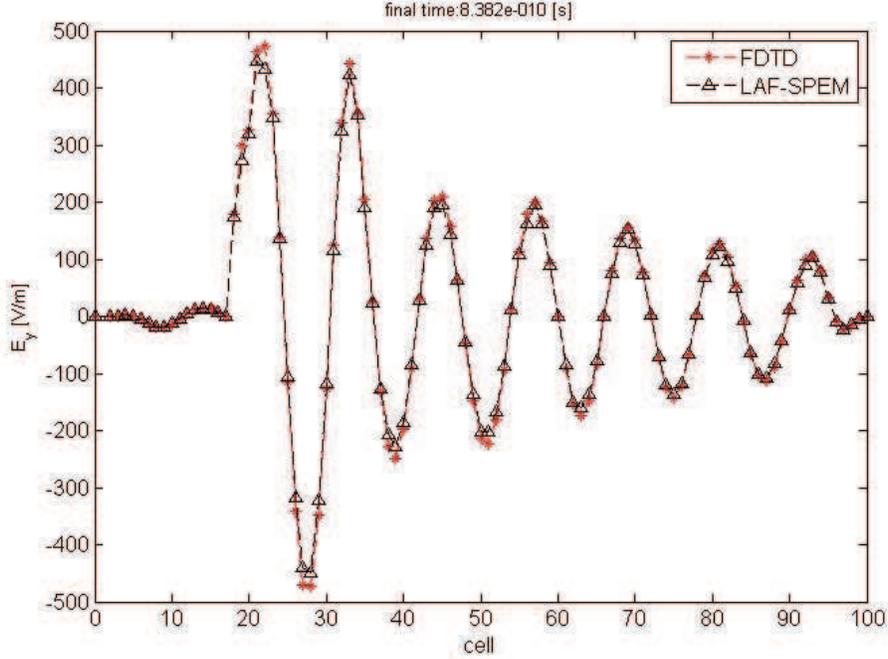}
\caption{Sectorial (2-D) PEC horn antenna excited by a sinusoidal voltage in a 2-D domain with 100x100 regular particles distribution. Comparison between FDTD and LAF-SPEM E field (V/m) results at time step 198, along the direction
containing the longitudinal axis of the antenna - relation (\ref{eq33})holds. }
\end{figure}

In order to show the capability of the proposed approach to run with a time step greater than that obtained by satisfying the CFL condition, a simulation with a time step equal to $4\Delta t$ CFL is carried out.
Also in this case, an irregular particles distribution has been also adopted by randomly positioning the particles near the original regular position in the domain.
As shown in figure 10, the simulation remains stable also for a final time corresponding to 200 time steps, and a satisfactory approximation has been achieved.

\begin{figure}[ht]
\label{fig:fig10}
\center
\includegraphics[width=12cm]{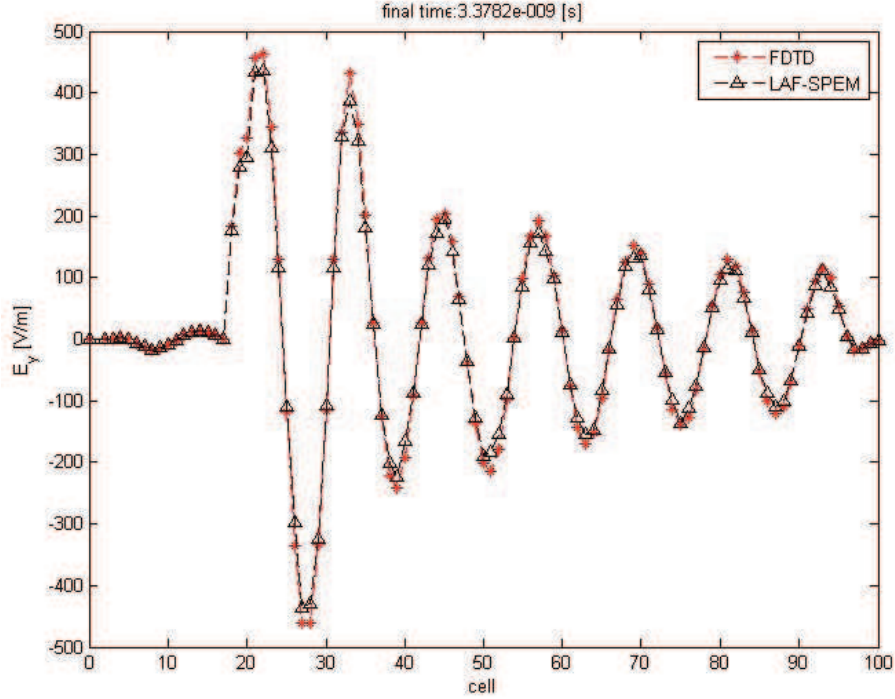}
\caption{ Sectorial (2-D) PEC horn antenna excited by a sinusoidal voltage in a 2-D domain with 100x100 irregular particles distribution. Comparison between $\textbf{E}_{y}$ field (V/m) results obtained with FDTD and LAF-SPEM at the same final time, along the direction containing the longitudinal axis of the antenna. The LAF-SPEM time step is 4 times that satisfying the CFL condition.}
\end{figure}

\section{Conclusions}
In this paper the SPH meshless method has been composed with a leapfrog ADI-FDTD method for
time evolution of electromagnetic time-domain transient propagation problems. The new schemeimproves the original SPEM method, already developed by the authors. The method is very promising
because avoids grid generation in space and results unconditionally stable in time. In this way, realistic simulations can be performed with reasonable computational efforts. The computational tool is assessed by using different numerical simulations also employing irregular particles distribution.  Perfect matching layer technique is used in simulating open spatial problems and a consistency restoring approach is introduced. The proposed methods may have good perspectives of extensive applications.\\

\end{document}